\documentclass[12pt,reqno]{amsart}
\usepackage{amssymb}
\usepackage{hyperref}

\allowdisplaybreaks

\newtheorem{thm}{Theorem}
\newtheorem{lem}{Lemma}
\newtheorem{cor}{Corollary}

\theoremstyle{definition}
\newtheorem{utv}{Proposition}
\newtheorem*{dfn}{Definition}

\theoremstyle{remark}
\newtheorem{rem}{Remark}

\numberwithin{equation}{section}

\newcommand{\vG}{\varGamma}
\newcommand{\veps}{\varepsilon}
\newcommand{\vphi}{\varphi}
\newcommand{\pwp}{\wp^{\prime}}
\newcommand{\wta}{\widetilde{a}}
\newcommand{\wtc}{\widetilde{c}}
\newcommand{\wtd}{\widetilde{d}}
\newcommand{\wtm}{\widetilde{m}}

\newcommand{\wtl}{\widetilde{\lambda}}
\newcommand{\wtI}{\widetilde{I}}

\newcommand{\wtG}{\widetilde{\vG}}
\newcommand{\whPsi}{\widehat{\Psi}}

\newcommand{\Bn}{\mathbb{N}}
\newcommand{\Bz}{\mathbb{Z}}

\newcommand{\Bc}{\mathbb{C}}
\newcommand{\ii}{\mathrm{i}}
\newcommand{\const}{\mathrm{const}}
\newcommand{\Ll}{\mathcal{L}}
\newcommand{\Cc}{\mathcal{C}}
\newcommand{\Pp}{\mathcal{P}}
\DeclareMathOperator{\im}{Im}

\newcommand{\lra}{\longrightarrow}
\def\nodot#1#2{}

\begin{document}
\title{Elliptic solitons and Heun's equation}
\author{A.~O. Smirnov}
\address{St.~Petersburg State University of Aircraft Engineering}
\email{smirnov@as.stoic.spb.su}
\translator{V.~B. Kuznetsov}


\begin{abstract}
We find a new class of algebraic geometric solutions of Heun's
equation with the accessory parameter belonging to a hyperelliptic
curve. Dependence of these solutions from the accessory parameter
as well as their relation to Heun's polynomials is studied.
Methods of calculating the algebraic genus of the curve, and its
branching points, are suggested. Monodromy group is considered.
Numerous examples are given.
\end{abstract}
\maketitle


\section*{Introduction}

The Heun equation (see, for instance, \cite{ODE,VTF, HDE, SL})
\begin{equation}\label{eq:heun}
\frac{d^2 y}{d z^2}+\left(\frac{\gamma}z+\frac{\delta}{z-1}
+\frac{\veps}{z-a}\right)\frac{d y}{d z} +\frac{\alpha\beta
z-q}{z(z-1)(z-a)}y=0,
\end{equation}
where
\begin{equation*}
1+\alpha+\beta-\gamma-\delta-\veps=0,
\end{equation*}
is a Fuchsian equation with four regular singularities at $0$, $1$, $a$
and $\infty$.  Recently, it started to attract much attention due to a
large number of its applications in mathematics and physics
\cite{KPS,Chn,AK,ZLS,KS,SLS}.

As a matter of fact, this equation has been known for a very long
time. Already in 1882 Darboux mentioned that it is the next one
beyond the hypergeometric equation and that ``it has already been
studied by many geometers'' \cite{D882}. In the works of Darboux
this equation appears when he investigated the minimal surfaces in
\cite{Dbk1} and when he studied harmonic equations in \cite{Dbk2}.
In the work \cite{D882} Darboux points out at a close connection
between the equation \eqref{eq:heun} and the generalization of
Lam\'e's equation
\begin{equation}
\psi_{xx}-u(x)\psi=E\psi, \label{eq:shr}
\end{equation}
\begin{gather}
u(x)=m_0(m_0+1)\wp(x)+\sum_{i=1}^3m_i (m_i+1)\wp(x-\omega_i),
\label{pot:tv}\\
\wp(\omega_i)=e_i, \qquad \wp(x-2\omega_i)\equiv \wp(x). \notag
\end{gather}
Here $\wp(x)$ is the Weierstrass function \cite{AKH},
\begin{equation}
[\pwp(x)]^2=4\wp^3(x)-g_2\wp(x)-g_3=4\prod_{j=1}^3(\wp(x)-e_j).
\label{eq:torus}
\end{equation}
In the same work Darboux remarks that the equation
\eqref{eq:shr}, \eqref{pot:tv} appeared already in the works
of Hermite on solution of a special case
(integer $N$) of Lam\'e's equation:
\begin{equation}
\psi_{xx}-N(N+1)\wp(x)\psi=E\psi. \label{eq:lame}
\end{equation}
Using the methods ``proposed by Hermite in 1872 in the course
of lectures in \'Ecole Polytechnique'', Darboux solves the equation
\eqref{eq:shr}, \eqref{pot:tv} with integer $m_i$ and points out that
one can find solutions of the equation \eqref{eq:shr}, \eqref{pot:tv}
for any value of $E$ in other special cases. In particular, ``when
all $m_i$ are half-odd-integer and $N$ is odd integer''
\footnote{so far we do not know by which method one can
find solutions of this special case of Heun's equation}.
Unfortunately, until recently the work \cite{D882} was not
well-known\footnote{the author learned about the work \cite{D882}
in 1999 from prof. V.B.Matveev}, in connection with this the
equation \eqref{eq:shr}, \eqref{pot:tv} has the name of Treibich
and Verdier who proved the finite-gapness of the potential
\eqref{pot:tv} for any integer $m_i$ \cite{TV90,TV90a,TV92}.

In the present paper we study exact solutions of a special case of
Heun's equation:
\begin{equation}\label{eq:heun1}
\frac{d^2 y}{d
z^2}+\frac{1}{2}\left(\frac{1-2m_1}z+\frac{1-2m_2}{z-1}
+\frac{1-2m_3}{z-a}\right)\frac{d y}{d z}
+\frac{N(N-2m_0-1)z+\lambda}{4z(z-1)(z-a)}y=0,
\end{equation}
where
\begin{gather}
N=m_0+m_1+m_2+m_3, \label{eq:N}\\
m_0,m_1,m_2,m_ 3 \in\Bz,\quad \lambda, z\in\Bc. \notag
\end{gather}

It is well-known (see, e.g., \cite{ODE, VTF, HDE, SL}) that by putting
\begin{equation*}
m_1=m_2=m_3=0, \qquad m_0=N
\end{equation*}
and changing the independent variable
\begin{subequations}\label{eq:z2x}
\begin{align}
&z=\frac{\wp(x)-e_1}{e_2-e_1},
&&a=\frac{e_3-e_1}{e_2-e_1}, \\
&y(z)=\psi(x),
&&\lambda=\frac{E-N(N+1)e_1}{e_1-e_2},
\end{align}
\end{subequations}
this equation is reduced to the familiar Lam\'e's equation
\eqref{eq:lame}. When $N\in\Bn$ the solutions of the equation
\eqref{eq:lame} are the eigenfunctions of the Schr\"o\-dinger
operator \eqref{eq:shr} with the $N$-gap elliptic potential
\begin{equation*}
u(x)=N(N+1)\wp(x).
\end{equation*}
They can be found either in terms of solutions of a differential
equation satisfied by products of these eigenfunctions
\cite{ODE,VTF,WW,D882} or by using the Krichever-Hermite ansatz
\cite{BBME,BBEIM,BE89b,BE94}.

In general, i.e. when $m_i\in\Bc$, the change of variables (see,
for instance, \cite{D882,Dbk1,Dbk2})
\begin{subequations}\label{eq:z2x1}
\begin{align}
&z=\frac{\wp(x)-e_1}{e_2-e_1},
&&a=\frac{e_3-e_1}{e_2-e_1}, \\
&y(z)=\psi(x)\prod_{i=1}^3(\wp(x)-e_i)^{m_i/2},
&&\lambda=\frac{E}{e_1-e_2}+\const,
\end{align}
\end{subequations}
transforms Heun's equation into Schr\"odinger's equation \eqref{eq:shr}
with the potential \eqref{pot:tv}, which for integer $m_i$ will be finite-gap
\footnote{changes of variables \eqref{eq:z2x},\eqref{eq:z2x1} are easily
generalized to the case of Schr\"odinger's equation with any even
elliptic potential. Based on our results in \cite{Sm89a,Smr94} we plan
a further study of a special case of the Fuchsian equation with five
singularities}.

Because the Treibich-Verdier equation \eqref{eq:shr},
\eqref{pot:tv} with integer $m_i$ has been actively studied over
the years by many researchers \cite{BE89b,BE94,TV90,TV90a,TV92},
including the author \cite{Sm89a,Sm89b,Smr94}, it became possible
to apply the obtained results to the theory of solutions of Heun's
equation.

In particular, in this paper it is proved that for any nonnegative
integer $m_i$ and for all $\lambda\ne\lambda_j$
($j=1,\ldots,2g+1$) the solutions of the Heun equation
\eqref{eq:heun1} are the functions of the following form:
\begin{equation*}
Y_{1,2}(m_0,m_1,m_2,m_3;\lambda;z)
=\sqrt{\Psi_{g,N}(\lambda,z)}\exp\left(\pm\frac{\ii\nu(\lambda)}2
\int\frac{z^{m_1}(z-1)^{m_2}(z-a)^{m_3}\,d
z}{\Psi_{g,N}(\lambda,z) \sqrt{z(z-1)(z-a)}}\right).
\end{equation*}
Here $\ii^2=-1$,
\begin{equation*}
\vG:\quad \nu^2=\prod_{j=1}^{2g+1}(\lambda-\lambda_j), \qquad
\lambda_j=\lambda(E_j),
\end{equation*}
$E_j$ are the gap edges of the finite-gap elliptic potential
$u(x)$ \eqref{pot:tv} and $\Psi_{g,N}(\lambda,z)$ is some
polynomial of the degree $N$ in $z$ and of the degree $g$ in
$\lambda$. We develop the methods of calculating the polynomial
$\Psi_{g,N}(\lambda,z)$, the branching points $\lambda_j$ of the
algebraic curve $\vG=\{(\nu,\lambda)\}$ and, respectively, the gap
edges $E_j$ of the spectrum of the Treibich-Verdier potential. We
find the dependence of the genus $g$ of the curves $\vG$ and
$\wtG=\{(w,E)\}$,
\begin{equation} \label{curve:hyp1}
\wtG:\quad w^2=\prod_{j=1}^{2g+1}(E-E_j),
\end{equation}
as a function of the characteristics (exponents)
$m_i$. The monodromy group of these solutions
is also considered.

We also prove that for any integer $m_i$'s and for
$\lambda=\lambda_j$ the solutions of Heun's equation are expressed
in terms of Heun polynomials \cite{VTF, HDE, SL}, i.e. Heun polynomials
with half-odd-integer $\beta-\alpha$, $\gamma$, $\delta$, $\epsilon$
are special cases
of the found solutions.

The needed facts from the theory of finite-gap elliptic potentials
for the Schr\"o\-dinger operator are collected without proof in the first
paragraph.

Because the exact solutions of the Heun equation are of much interest,
in the Appendix we give several simplest solutions of the Heun equation
belonging to this class.

The author thanks S.Yu.~Slavyanov, V.B.Matveev, V.B.~Kuznetsov and
A.Ya.~Ka\-zakov for useful discussions and the referees for
helpful suggestions about improving the text.

\section{Schr\"odinger operator with finite-gap elliptic potential}

\begin{utv}[\cite{GW94d, TV97, TV99}] \label{utv:pot:tv}
For any $m_i \in \Bz$ the Treibich-Verdier potential
\eqref{pot:tv} is a finite-gap elliptic potential for the
Schr\"odinger operator \eqref{eq:shr}. For any nonnegative integer
$m_i$ the following formula for the genus $g$ of the spectral
curve $\wtG$ \eqref{curve:hyp1} for the potential \eqref{pot:tv}
is true:
\begin{equation}
g=\frac{1}{2}\max\left\{2\max_{0\leqslant i\leqslant 3} m_i,
1+N-(1+(-1)^N)\left(\min_{0\leqslant i\leqslant 3} m_i
+\frac{1}{2}\right)\right\}.
\end{equation}
\end{utv}

\begin{utv}[\cite{Nov74,ZMNP,KD,Mat76}]
Any $g$-gap potential $u(x)$ of the Schr\"odinger operator
\eqref{eq:shr} satisfies the Novikov equation:
\begin{equation}
J_g+\sum_{m=1}^g c_m J_{g-m}=d, \label{eq:novikov}
\end{equation}
or, which is the same, the stationary `higher' Korteweg-de Vries
(KdV) equation:
\begin{equation*}
\partial_x\left(J_g+\sum_{m=1}^g c_m J_{g-m}\right)=0.
\end{equation*}
Here $c_m$, $d$ are constants and the functions $J_m$
are the flows of the `higher' KdV equations
\begin{equation*}
\partial_{t_m}u=\partial_x J_m.
\end{equation*}
The expressions for the flows $J_m$ are found from the
relations
\begin{subequations} \label{kdv:int:J}
\begin{align}
&L\psi=\psi_{xx}-4u\psi+2u_x\int_{x}^{\infty}\psi(\tau)d\tau, \\
&(J_n)_x=L^n(u_x).
\end{align}
\end{subequations}
where $u(x)$ is a potential decreasing fast at $\infty$. In particular,
\begin{align*}
&J_0=u, \qquad  J_1=u_{xx}-3u^2, \qquad
J_2=u_{xxxx}-10u_{xx}u-5u^2_x+10u^3, \\
&J_3=u_{xxxxxx}-14uu_{xxxx}-28u_xu_{xxx}-21u_{xx}^2
+70u^2u_{xx}+70uu^2_x- 35u^4.
\end{align*}
In the case of decreasing fast at $\infty$ potential the variables
\begin{equation*}
C_j=\int_{-\infty}^{\infty} J_j(x,t)\,d x
\end{equation*}
constitute an infinite set of integrals of motion for the KdV equation
\begin{equation*}
\partial_t u=\partial_x J_1.
\end{equation*}
\end{utv}

\begin{utv}[\cite{DMN76,Mat76,ZMNP}]
The function
\begin{equation}
\whPsi(x,E)=E^g+\sum_{j=1}^g \gamma_j(x) E^{g-j},
\label{kdv:prod:psi}
\end{equation}
where
\begin{equation}
\gamma_j(x)=-\frac2{4^j}\left(J_{j-1}
+\sum_{m=1}^{j-1}c_m J_{j-m-1}-\frac{c_j}2\right), \label{kdv:gamma}
\end{equation}
obeys the equation
\begin{equation*}
\whPsi_{xxx}=4(u+E)\whPsi_x+2u_x\whPsi,
\end{equation*}
the solutions of which are the products of any two solutions of the
equation \eqref{eq:shr}. Here $c_j$, $J_j$ are the same as in
\eqref{eq:novikov} and $u(x)$ is a $g$-gap potential of the
Schr\"odinger operator \eqref{eq:shr}.
\end{utv}

\begin{utv}[\cite{BBME,BBEIM,BE89b,BE94}]
The equation \eqref{eq:shr} with the potential of
Treibich-Verdier \eqref{pot:tv} has a solution
of the form (Krichever-Hermite ansatz)
\begin{equation}
\psi(x;k,\alpha)=\sum_{i=0}^3\sum_{j=0}^{m_i-1}
a_{ij}(k,\alpha)e^{kx}\frac{d^j}{dx^j}\Phi(x-\omega_i,\alpha).
\label{psi:sum}
\end{equation}
Here $k$, $\alpha$ are auxiliary spectral parameters,
$\omega_0\equiv 0$,
\begin{equation}
\Phi(x,\alpha)=\frac{\sigma(\alpha-x)}{\sigma(x)\sigma(\alpha)}
\exp\{\zeta(\alpha)x\},
\label{eq:Phi}
\end{equation}
and $\sigma(\alpha)$, $\zeta(\alpha)$ are Weierstrass' elliptic
functions \cite{AKH}.

Substituting the ansatz \eqref{psi:sum} for the $\psi$-function into
the equation \eqref{eq:shr} one can find the equation of the spectral
curve $\vG$ for the finite-gap potential $u(x)$
\begin{equation*}
\vG=\{(k,\alpha)\}:\qquad R(k,\alpha)=0,
\end{equation*}
as well as the dependence of the coefficients $a_{ij}(k,\alpha)$
and the main spectral parameter $E$
\begin{equation*}
E(-k,-\alpha)=E(k,\alpha)
\end{equation*}
as functions of the auxiliary spectral parameters $k$ and
$\alpha$.

The gap edges of the spectrum of the finite-gap potential $u(x)$
are the values of the spectral parameter $E$ calculated at the
points $\Pp$ of the curve $\vG$
\begin{equation*}
E_j=E(\Pp_j),\qquad \tau \Pp_j=\Pp_j,\qquad j=1,\ldots,2g+2,
\qquad E_{2g+2}=\infty.
\end{equation*}
which are invariant with respect to the involution $\tau$
\begin{equation*}
\tau:\;(k,\alpha)\lra(-k,-\alpha).
\end{equation*}

Linear independent solutions of the equation \eqref{eq:shr} for all values
of the spectral parameter $E$, except for the gap edges $E_j$, are the
functions
\begin{equation}
\psi_1(x,E)=\psi(x;k,\alpha)\quad\text{and}\quad
\psi_2(x,E)=\psi(x;-k,-\alpha).
\label{eq:psi12}
\end{equation}

A product of these functions is a doubly periodic function and, modulo
the normalizing factor, it is equal to $\whPsi(x,E)$ \eqref{kdv:prod:psi},
\eqref{kdv:gamma}.
\end{utv}

\section{Finite-gap solutions of Heun's equation}

Since the potential \eqref{pot:tv} is finite-gap, one can introduce
the notion of {\it `finite-gap' solutions} to the Heun equation which
are obtained from the finite-gap solutions of the Treibich-Verdier
equation by the inverse change of variables  (see \eqref{eq:z2x})
\begin{subequations} \label{eq:x2z}
\begin{gather}
\wp(x)=e_1+(e_2-e_1)z,\qquad
e_2=\frac{a-2}{a+1}e_1,\qquad
e_3=\frac{1-2a}{a+1}e_1,\\
E=(e_1-e_2)\lambda+\const,\qquad
u(x)=4(e_2-e_1)U(z)+\const_1. \label{eq:x2z2}
\end{gather}
\end{subequations}

\begin{dfn} By {\it Novikov's equation for the Heun equation
of order $g$} we will call the following equality:
\begin{equation}
I_g+\sum_{j=1}^g \wtc_j I_{g-j}=\wtd, \label{eq:nov1}
\end{equation}
where $\wtc_j$, $\wtd$ are some constants,
\begin{subequations} \label{heun:int:I}
\begin{equation}
I_0=U(z), \qquad I_{j+1}=\Ll(I_j),
\end{equation}
\begin{equation}
\Ll(f)=z(z-1)(z-a)\frac{d^2f}{dz^2}
+\frac{3z^2-2(a+1)z+a}2\cdot\frac{d f}{d z}
-\int\left(4I_0\frac{d f}{d z}+2f\frac{dI_0}{d z}\right)d z,
\end{equation}
\end{subequations}
$U(z)$ is the `potential' \eqref{eq:x2z2}\footnote{the term
``potential'' attached to the function $U(z)$ is used here
only in analogy with the function $u(x)$, and it does not carry
any additional meaning}
\begin{multline} \label{pot:heun}
U(z)=\frac{m_0(m_0+1)}4z +\frac{m_1(m_1+1)}4\cdot\frac a z+{}\\
{}+\frac{m_2(m_2+1)}4\cdot\frac{z-a}{z-1}
+\frac{m_3(m_3+1)}4\cdot\frac{a(z-1)}{z-a}.
\end{multline}
\end{dfn}

\begin{thm}
The `potential' $U(z)$ \eqref{pot:heun} satisfies Novikov's equation
for the Heun equation of order $g$ \eqref{eq:nov1}, \eqref{heun:int:I}
if and only if the potential $u(x)$  \eqref{pot:tv} satisfies Novikov's
equation \eqref{eq:novikov}, \eqref{kdv:int:J} of the same order.
\end{thm}

\begin{dfn}  `Potential' $U(z)$ will be called $g$-gap potential if it satisfies
Novikov's equation for the Heun equation of order $g$
\eqref{eq:nov1}, \eqref{heun:int:I}.
\end{dfn}

\begin{rem} The flows $J_m$ \eqref{kdv:int:J} in the Novikov equation
\eqref{eq:novikov} were originally defined for the potentials decreasing
fast at $\infty$. Although for finite-gap as well as for any periodic
potentials these functions can not be obtained from the equations
\eqref{kdv:int:J} because of the divergent integral with an infinite upper
limit, the expressions for $J_m$ are still valid in the case of
a finite-gap potential.
The constant of integration in the definition of the `flows' $I_n$
\eqref{heun:int:I} can not be easily fixed too. Because of that,
all the `flows' $I_n$ are defined modulo a linear combination
of lower order `flows'. However, it is easily seen that the property
of `finite-gapness' of the `potential' $U(z)$ \eqref{pot:heun}
does not depend from the concrete values of integration constants
in the definition of $I_n$. The latter affect only the values of the
constants $\wtc_m$ and $\wtd$ in the equation \eqref{eq:nov1}.
\end{rem}

\begin{thm} For any $m_i\in\Bz$ there exists a nonnegative
integer $g$ such that the `potential' $U(z)$ \eqref{pot:heun}
satisfies Novikov's equation for the Heun equation of order $g$
\eqref{eq:nov1}, \eqref{heun:int:I}.
\end{thm}

\begin{proof} Assume that all characteristics $m_i$ in the equation
\eqref{pot:heun} are nonnegative. If anyone $m_i<0$ then one can
always make the transformation $m_i \to -m_i-1$ which does not change
the `flows' $I_n$ \eqref{heun:int:I}, \eqref{pot:heun}.

From the properties of the flows $J_n$ \eqref{kdv:int:J} (see,
e.g., \cite{KD}), from the properties of elliptic functions
\cite{AKH} and from the equation of the change of variable
\eqref{eq:x2z} it follows that all `flows' $I_n$
\eqref{heun:int:I} are rational functions of the variable $z$
(i.e. these functions do not have logarithmic singularities).

It is also not very difficult to check that the order of poles of
the functions $I_n(z)$ in the singularities of Heun's equation
does not exceed the corresponding characteristic $m_i$.

Indeed, it is easy to verify that for $m_1=0$ neither
$I_0(z)=U(z)$ \eqref{pot:heun} nor all other $I_n(z)$
have a pole at $z=0$, but for $m_1 \neq 0$ the function
$I_0(z)$ has at this point a pole of first order.

If we now assume that the function $I_n(z)$ has at the point $z=0$
a pole of order $\alpha\leqslant m_1$ (which is undoubtedly true
for $n=0$ and $m_1\neq 0$):
\begin{equation*}
I_n(z)=\frac A{z^{\alpha}}+O\left(z^{1-\alpha}\right),
\qquad  z\to 0,
\end{equation*}
then we obtain that the function $I_{n+1}(z)$,
\begin{equation*}
I_{n+1}(z)=\Ll\big(I_n(z)\big)=
\frac{(2\alpha+1)(\alpha-m_1)(\alpha+m_1+1)}{2(\alpha+1)} \cdot
\frac{a A}{z^{\alpha+1}} +O\left(z^{-\alpha}\right),
\end{equation*}
has at the same point a pole of order
$\alpha'=\min\{m_1,\alpha+1\}$. By similar means one can check the
other poles of the function $I_n(z)$.

Hence, for any $n$ the dimension of the linear span
of rational functions $1,I_0,\ldots,I_n$ does not exceed $N+1$ and therefore
there exists a number $g$,
\begin{equation}
\max_{0\leqslant i\leqslant 3} m_i \leqslant g\leqslant N,
\label{eq:g}
\end{equation}
such that the equality \eqref{eq:nov1} is fulfilled.
\end{proof}

\begin{cor} \label{cor:pot:tv}
For any $m_i\in\Bz$ the potentials $U(z)$ \eqref{pot:heun} and $u(x)$
\eqref{pot:tv} are finite-gap. The genus $g$ of the spectral curve
$\vG$, associated with the finite-gap potential $u(x)$
\eqref{pot:tv}, satisfies the inequality \eqref{eq:g}.
\end{cor}

\begin{dfn} Solutions of the Heun equation \eqref{eq:heun1}
with integer characteristics $m_i$ will be called {\it finite-gap solutions
of Heun's equation}.
\end{dfn}

\begin{rem} As it was said above (cf. Proposition \ref{utv:pot:tv})
the proof of finite-gapness of the potential $u(x)$ \eqref{pot:tv}
was first given in the works by Treibich and Verdier \cite{TV90a,TV92}.
Here we gave another proof which allowed us to specify the
explicit form of the finite-gap solutions of Heun's equation.
\end{rem}

In order to find solutions of Heun's equation \eqref{eq:heun1}
let us use the results of the theory of
finite-gap elliptic potentials for the Schr\"odinger operator and
consider the equation ($m_i\geqslant 0$, $i=0,1,2,3$)
\begin{equation}
\frac{d^3\Psi}{dz^3}+3P(z)\frac{d^2\Psi}{dz^2}
+\left(P'(z)+4Q(z)+2P^2(z)\right)\frac{d\Psi}{d z}
+\left(2Q'(z)+4P(z)Q(z)\right)\Psi=0,\label{eq:3}
\end{equation}
\begin{align*}
P(z)&=\frac12\left(\frac{1-2m_1}z+\frac{1-2m_2}{z-1}
+\frac{1-2m_3}{z-a}\right),\\
Q(z)&=\frac{N(N-2m_0-1)z+\lambda}{4z(z-1)(z-a)},
\end{align*}
solutions of which are the products of any two solutions of
Heun's equation \eqref{eq:heun1}.

\begin{thm} \label{thm:heun:Psi} Equation \eqref{eq:3} with nonnegative integer
characteristics $m_i$ has as its solution the function $\Psi_{g,N}(\lambda,z)$,
which is a polynomial in $\lambda$ of the degree $g$ and in $z$
of the degree $N$ \eqref{eq:N}. The leading coefficient of this function
considered as a polynomial in $\lambda$ is equal
\begin{equation}
\wta_0(z)=z^{m_1}(z-1)^{m_2}(z-a)^{m_3}.
\label{heun:a0}
\end{equation}
In other words,
\begin{equation} \label{eq:heun:Psi}
\begin{split}
\Psi_{g,N}(\lambda,z)&=a_0(\lambda) z^N+a_1(\lambda) z^{N-1}
+\ldots+a_N(\lambda)={}\\
&{}=z^{m_1}(z-1)^{m_2}(z-a)^{m_3}\lambda^g
+\wta_1(z)\lambda^{g-1}+\ldots+\wta_g(z).
\end{split}
\end{equation}
\end{thm}

\begin{proof} From \eqref{kdv:prod:psi}--\eqref{eq:psi12} it follows that the
product $\whPsi(x,E)$ of eigenfunctions of the Schr\"odinger operator
with the Treibich-Verdier potential is an elliptic meromorphic function
in the variable $x$, and it is a polynomial of degree $g$ in the spectral
parameter $E$. As a function of the variable $x$ the function
$\whPsi(x,E)$ has poles of multiplicity $2m_i$ at the points
$x=\omega_j$ ($\omega_0\equiv 0$).

From \eqref{kdv:prod:psi}, \eqref{pot:tv} and from the properties
of the Weierstrass $\wp$-function it follows that $\whPsi(x,E)$ is
a rational function in $\wp(x)$. Hence, the function
\begin{equation}
\Psi_{g,N}(\lambda,z)=\const\cdot\whPsi(x,E)
\prod_{j=1}^3(\wp(x)-e_j)^{m_i},
\label{heun:Psi:def}
\end{equation}
where $\lambda$ and $z$ are related with $E$ and $x$
by equalities \eqref{eq:x2z}, is a polynomial in $\lambda$ of degree
$g$ and in $z$ of degree $N$ (i.e. it is a rational function in the variable $z$
with the unique pole of order $N$ at the point $z=\infty$). The constant
in the equality \eqref{heun:Psi:def} is chosen such that the leading
coefficient of $\Psi_{g,N}(\lambda,z)$, considered as a polynomial
in $\lambda$, to be equal \eqref{heun:a0}.
\end{proof}

\begin{cor}
Coefficients $\wta_j(z)$ of the polynomial $\Psi_{g,N}(\lambda,z)$
have the form:
\begin{equation}
\wta_j(z)=z^{m_1}(z-1)^{m_2}(z-a)^{m_3}\wtI_{j-1},
\qquad j=1,\ldots,g,
\end{equation}
where $\wtI_{j}$ is a linear combination of the rational functions
$I_j,\ldots,I_0,1$ having poles in the singularities of Heun's
equation.
\end{cor}

\begin{proof}[Proof\nodot] follows from the
equalities \eqref{kdv:gamma}, \eqref{heun:Psi:def} and from
the change of variable \eqref{eq:x2z}.
\end{proof}

Knowing the product of the solutions of Heun's equation
it is not difficult to find the solutions themselves (see, for instance,
\cite[\S 19.53, \S 23.7, \S23.71]{WW}).

\begin{thm} Finite-gap solutions of Heun's equation with nonnegative
characteristics $m_i$ have the form
\begin{multline} \label{eq:Y}
Y_{1,2}(m_0,m_1,m_2,m_3;\lambda;z)={}\\
{}=\sqrt{\Psi_{g,N}(\lambda,z)}\exp\left(\pm\frac{\ii\nu(\lambda)}2
\int\frac{z^{m_1}(z-1)^{m_2}(z-a)^{m_3}\,d
z}{\Psi_{g,N}(\lambda,z) \sqrt{z(z-1)(z-a)}}\right).
\end{multline}
Here $\ii^2=-1$,
\begin{equation}
\vG:\quad \nu^2=\prod_{j=1}^{2g+1}(\lambda-\lambda_j), \qquad
\lambda_j=\lambda(E_j), \label{curve:hyp}
\end{equation}
$E_j$ are the gap edges of the finite-gap elliptic potential
$u(x)$ \eqref{pot:tv}.
\end{thm}

\begin{proof}
From the Liouville formula it follows that the Wronskian of two
linearly independent solutions of the linear homogeneous
differential equation,
\begin{equation*}
y''+P(z)y'+Q(z)y=0,
\end{equation*}
has the following dependence from $z$:
\begin{equation*}
W[y_1(\lambda,z),y_2(\lambda,z)]=W_0(\lambda)
\exp\left\{-\int_{z_0}^z P(t)\,dt\right\},
\end{equation*}
where $W_0(\lambda)$ is Wronskian's value at $z_0$.
Hence, the Wronskian of two linearly independent solutions
of Heun's equation \eqref{eq:heun1} is equal to
\begin{equation}\label{eq:wron}
W[y_1(\lambda,z),y_2(\lambda,z)]=-\ii\nu(\lambda)\cdot
\frac{z^{m_1}(z-1)^{m_2}(z-a)^{m_3}}{\sqrt{z(z-1)(z-a)}},
\end{equation}
where
\begin{equation*}
\nu(\lambda)= \ii W_0(\lambda)\cdot\frac{\sqrt{z_0(z_0-1)(z_0-a)}
}{z_0^{m_1}(z_0-1)^{m_2}(z_0-a)^{m_3}}.
\end{equation*}
If we now divide the Wronskian of two solutions \eqref{eq:wron}
by their product,
\begin{equation}\label{eq:prod}
y_1(\lambda,z)\cdot y_2(\lambda,z)=\Psi_{g,N}(\lambda,z),
\end{equation}
we obtain a simple differential equation of the first oder:
\begin{equation}\label{eq:frac:dif}
\frac{y_2'}{y_2}-\frac{y_1'}{y_1}= -\ii\nu(\lambda)\cdot
\frac{z^{m_1}(z-1)^{m_2}(z-a)^{m_3}}{\Psi_{g,N}(\lambda,z)\sqrt{z(z-1)(z-a)}}.
\end{equation}
The equation \eqref{eq:frac:dif} can be easily integrated:
\begin{equation}\label{eq:frac}
\frac{y_2(\lambda,z)}{y_1(\lambda,z)}
=C(\lambda)\cdot\exp\left[-\ii\nu(\lambda)\int_{z_1}^z
\frac{t^{m_1}(t-1)^{m_2}(t-a)^{m_3}\,dt}{\Psi_{g,N}(\lambda,t)
\sqrt{t(t-1)(t-a)}}\right],
\end{equation}
where $C(\lambda)=y_2(\lambda,z_1)/y_1(\lambda,z_1)$.

If we now consider solutions $Y_{1,2}(\lambda,z)$ with the same
product \eqref{eq:prod} and Wronskian \eqref{eq:wron}
\begin{equation}
Y_1(\lambda,z)=\sqrt{C(\lambda)}\cdot y_1(\lambda,z),\qquad
Y_2(\lambda,z)= \frac{y_2(\lambda,z)}{\sqrt{C(\lambda)}},
\end{equation}
then from \eqref{eq:prod} and \eqref{eq:frac} we get \eqref{eq:Y}.

Substituting the ansatz \eqref{eq:Y} into the Heun equation
\eqref{eq:heun1} we get
\begin{equation}
\nu^2(\lambda)=\frac{2\Psi\Psi''-(\Psi')^2+2P(z)\Psi\Psi'
+4Q(z)\Psi^2}{z^{2m_1-1}(z-1)^{2m_2-1}(z-a)^{2m_3-1}},
\label{eq:heun:nu}
\end{equation}
where
\begin{equation*}
\Psi=\Psi_{g,N}(\lambda,z),\qquad \Psi'= \frac{d\Psi}{d z},\qquad
\Psi''= \frac{d^2\Psi}{d z^2}.
\end{equation*}
From \eqref{eq:heun:Psi}, \eqref{eq:heun:nu} it follows that
$\nu^2(\lambda)$ is a polynomial in $\lambda$ of the degree $2g+1$
with the leading coefficient equal to~$1$.

It is not difficult to show that under the change  \eqref{eq:z2x1}
the solutions $Y_{1,2}(\lambda,z)$
of Heun's equation \eqref{eq:heun1} turn into Floquet solutions
of the Treibich-Verdier equation \eqref{eq:shr},
\eqref{pot:tv}. Therefore, zeros $\lambda_j$ ($j=1,\ldots,2g+1$)
of the polynomial $\nu(\lambda)$ or, which is the same, zeros
of the Wronskian of the solutions
$Y_{1,2}(\lambda,z)$ correspond to zeros of the Wronskian of
the Floquet solutions of the Treibich-Veridier equation, i.e. they correspond
to the
gap edges $E_j$ ($j=1,\ldots,2g+1$) of spectrum of the Triebich-Verdier
potential. Hence, the hyperelliptic curve $\vG$ \eqref{curve:hyp} is
isomorphic to the spectral curve $\wtG$ \eqref{curve:hyp1} of the
finite-gap elliptic potential $u(x)$ \eqref{pot:tv}.
\end{proof}

\begin{rem} The fact that the equation \eqref{eq:3} for any nonnegative
integer $m_i$ and for any $\lambda$ has as its solution a polynomial
in $z$ of degree $N$ was known to Darboux who in the work \cite{D882}
pointed put at the method of finding some exact solutions of the generalized
Lam\'e equation, which is now known as the `Treibich-Verdier equation'.
However, he did not study the dependence of this polynomial and, respectively,
solutions of the Heun equation, from $\lambda$.
\end{rem}

\begin{rem}
In order to find the `finite-gap' solutions of Heun's equation
with negative characteristics one must use the equalities (see
\cite{ODE, VTF, HDE, SL}).
\begin{enumerate}
\item Negative $m_1$:
\begin{equation} \label{eq:-m1}
Y(m_0,-m_1-1,m_2,m_3;\lambda;z)=z^{-m_1-1/2}
Y(m_0,m_1,m_2,m_3;\mu;z),
\end{equation}
where
\begin{equation*}
\mu=\lambda-(2m_1+1)(2m_2-1)a-(2m_1+1)(2m_3-1).
\end{equation*}
\item Negative $m_2$:
\begin{equation}
Y(m_0,m_1,-m_2-1,m_3;\lambda;z)=(z-1)^{-m_2-1/2}
Y(m_0,m_1,m_2,m_3;\mu;z),
\end{equation}
where
\begin{equation*}
\mu=\lambda-(2m_1-1)(2m_2+1)a.
\end{equation*}
\item Negative $m_3$:
\begin{equation}
Y(m_0,m_1,m_2,-m_3-1;\lambda;z)=(z-a)^{-m_3-1/2}
Y(m_0,m_1,m_2,m_3;\mu;z),
\end{equation}
where
\begin{equation*}
\mu=\lambda-(2m_1-1)(2m_3+1).
\end{equation*}
\item Negative $m_0$:
\begin{equation} \label{eq:-m0}
Y(-m_0-1,m_1,m_2,m_3;\lambda;z)=Y(m_0,m_1,m_2,m_3;\lambda;z).
\end{equation}
\end{enumerate}

\end{rem}

\section{Finite-gap solutions and Heun polynomials}

Looking at the analytical properties of the solutions \eqref{eq:Y}
it is not difficult to prove the following statements.

\begin{cor}\label{cor:heun:zero1} For $\lambda\neq\lambda_j$
$(j=1,\ldots,2g+1)$, where $\lambda_j$ are the branching points
of the hyperelliptic curve $\vG$ \eqref{curve:hyp}\textup{:}
\begin{enumerate}
\item the polynomial $\Psi_{g,N}(\lambda,z)$ does not have zeros at the
singularities of Heun's equation;
\item the leading coefficient $a_0(\lambda)$ of the polynomial $\Psi_{g,N}(\lambda,z)$
is not equal to zero;
\item all zeros $z_k$ $(k=1,\ldots,N)$ of the polynomial
$\Psi_{g,N}(\lambda,z)$ are simple zeros satisfying the relations
\begin{equation*}
\prod_{\substack{j=1\\j\ne k}}^N(z_k-z_j)^2=
-\frac{\nu^2(\lambda)}{a_0^2(\lambda)}z_k^{2m_1-1}(z_k-1)^{2m_2-1}(z_k-a)^{2m_3-1}.
\end{equation*}
\end{enumerate}
\end{cor}

\begin{cor} Zeros of the equations
\begin{subequations}\label{eq:lj}
\begin{align}
&a_0(\lambda)=0,\\
&\Psi_{g,N}(\lambda,0)=0,\\
&\Psi_{g,N}(\lambda,1)=0,\\
&\Psi_{g,N}(\lambda,a)=0,\\
&\Delta(\lambda)=0,
\end{align}
\end{subequations}
where $\Delta(\lambda)$ is the discriminant of the algebraic equation
of $N$th order in $z$:
\begin{equation*}
\Psi_{g,N}(\lambda,z)=0,
\end{equation*}
are among the branching points $\lambda_j$ of the curve $\vG$
\eqref{curve:hyp}.
\end{cor}

\begin{lem}
For $\lambda=\lambda_j$ the solutions \eqref{eq:Y} cease to be
linearly independent. In this case the linear independent
`finite-gap' solutions of Heun's equation \eqref{eq:heun1} are
the functions
\begin{subequations} \label{eq:Y1}
\begin{align}
Y_1(m_0,m_1,m_2,&m_3;\lambda_j;z)=\sqrt{\Psi_{g,N}(\lambda_j,z)}\\
\intertext{and}
Y_2(m_0,m_1,m_2,&m_3;\lambda_j;z)=\sqrt{\Psi_{g,N}(\lambda_j,z)}
\int\frac{z^{m_1}(z-1)^{m_2}(z-a)^{m_3}\,d
z}{\Psi_{g,N}(\lambda_j,z) \sqrt{z(z-1)(z-a)}}.
\end{align}
\end{subequations}
\end{lem}

\begin{cor}\label{cor:lj} For $\lambda=\lambda_j$ $(j=1,\ldots,2g+1)$,
where $\lambda_j$ are the branching points of the hyperelliptic curve
$\vG$ \eqref{curve:hyp}\textup{:}
\begin{enumerate}
\item
the polynomial $\Psi_{g,N}(\lambda_j,z)$ has the form
\begin{equation}\label{eq:Psi.lj}
\Psi_{g,N}(\lambda_j,z)=z^{M_{1j}}(z-1)^{M_{2j}}
(z-a)^{M_{3j}}F^2_{nj}(z),
\end{equation}
where
$M_{ij}\in\{0,2m_i+1\}$ and $F_{nj}(z)$ is a polynomial in $z$
of order $n$;
\item the polynomial $F_{nj}(z)$ is a solution of Heun's equation with
the $P$-symbol
\begin{equation}
P\begin{Bmatrix}
0&1&a&\infty&;z\hfill\\
0&0&0&-n&;-\wtl_j/4\\
1/2+\wtm_{1j}&1/2+\wtm_{2j}&1/2+\wtm_{3j}&-n+1/2+\wtm_{0j}
\end{Bmatrix},
\end{equation}
where
\begin{equation}\label{eq:n}
\begin{split}
&M_{0j}\in\{0,2m_0+1\},\\
&\wtm_{ij}=m_i-M_{ij},\qquad
\wtm_{ij}\in\{m_i,-m_i-1\},\qquad i=0,1,2,3,\\
&n=(\wtm_{0j}+\wtm_{1j}+\wtm_{2j}+\wtm_{3j})/2,\\
&\wtl_j=\lambda_j+2m_1M_{3j}+2m_3M_{1j}-2M_{1j}M_{3j}-M_{1j}
-M_{3j}+{}\\
&{}\qquad+(2m_1M_{2j}+2m_2M_{1j}-2M_{1j}M_{2j}-M_{1j}-M_{2j})a,
\end{split}
\end{equation}
i.e. the polynomial $F_{nj}(z)$ modulo a constant factor is
a Heun polynomial \textup{(}see, e.g., \cite{HDE, SL}\textup{)};
\item the order of the polynomial $\Psi_{g,N}(\lambda_j,z)$, as a
polynomial in $z$, is equal to $(N-M_{0j})$;
\item
all zeros $z_k$ $(k=1,\ldots,n)$ of the polynomial $F_{nj}(z)$
are simple and satisfy the relations
\begin{equation*}
\sum_{\substack{i=1\\i\ne k}}^n\frac1{z_k-z_i}=\frac14
\left(\frac{2\wtm_{1j}-1}{z_k}+\frac{2\wtm_{2j}-1}{z_k-1}+
\frac{2\wtm_{3j}-1}{z_k-a}\right).
\end{equation*}
\end{enumerate}
\end{cor}

As a consequence, to every branching point $\lambda_j$
$(j=1,\ldots,2g+1)$ there corresponds a set of numbers $\wtm_{ij}$
$(i=0,1,2,3)$ such that $n$ \eqref{eq:n} is a nonnegative integer.
Solutions of Heun's equation for $\lambda=\lambda_j$ are expressed
in terms of Heun's polynomials of order $n$.

\begin{cor} For $N\neq 0$ every branching point $\lambda_j$
of the hyperelliptic curve $\vG$ \eqref{curve:hyp} is a zero of
one of the equations \eqref{eq:lj}.
\end{cor}

\begin{proof}[Proof\nodot] is done by assuming the opposite.
Let there exists a branching point $\lambda_j$ which is not a zero of anyone
of the equations \eqref{eq:lj}. From the condition it follows that in this
case the singularities of Heun's equation are not the zeros of the polynomial
$\Psi_{g,N}(\lambda_j,z)$ \eqref{eq:Psi.lj}. Therefore
\begin{equation}\label{eq:Mij=0}
M_{1j}=M_{2j}=M_{3j}=0.
\end{equation}

From the equalities \eqref{eq:Psi.lj}, \eqref{eq:Mij=0} and from
the condition $\Delta(\lambda_j)\neq0$ it follows that
$\Psi_{g,N}(\lambda_j,z)$ $= \const$. But this, for $N\neq 0$, is
in contradiction with the condition $a_0(\lambda_j)\neq0$. Hence,
every branching point $\lambda_j$ of the hyperelliptic curve $\vG$
\eqref{curve:hyp} for $N\neq 0$ is a zero of one of the equations
\eqref{eq:lj}.
\end{proof}

\begin{lem} \label{lem:nk}
For any set of integers $m_i$ $(i=0,1,2,3)$ there exists
not less than one and not more than four numbers
$n_k$ such that
\begin{equation} \label{eq:nk}
\begin{split}
n_k=1+\frac12\sum_{i=0}^3 \wtm_{ik},&\qquad n_k\in\Bn,\\
\wtm_{ik}\in\{m_i,-m_i-1\},&\qquad i=0,1,2,3.
\end{split}
\end{equation}
\end{lem}

\begin{proof} Rewrite the equality \eqref{eq:nk} in the form
\begin{equation}
\begin{split}
n_k&{}=\frac12\sum_{i=0}^3 \left(\wtm_{ik}+\frac12\right)={}\\
&{}=\frac12\sum_{i=0}^3(-1)^{\veps_{ik}}\left(m_i+\frac12\right),
\qquad \veps_{ik}\in\{0,1\}.
\end{split}
\end{equation}
It is easy to see that from sixteen values which the number
$n_k$ can take there are only eight which are integer, and that
to every positive number $n_k$ there corresponds a sign-opposite
negative one. Therefore the number of positive integer $n_k$'s can be
not more than four.

On the other hand, all integer $n_k$ can not be equal to zero because,
as it is not difficult to show, this was possible only if all
$m_i=-1/2$, $(i=0,1,2,3)$ which is in contradiction with the condition
$m_i\in\Bz$. Hence, among integer $n_k$ there is at least one positive
and one negative number.
\end{proof}

\begin{thm}\label{thm:polyheun}
Let $m_i$ $(i=0,1,2,3)$ be an arbitrary set of integers
satisfying the condition:
\begin{equation}\label{eq:mi}
n=1+\frac12\sum_{i=0}^3 m_i,\qquad n\in\Bn.
\end{equation}
Then:
\begin{enumerate}
\item there are exactly $n$ values $\lambda_j$ such that for
$\lambda=\lambda_j$ the Heun equation \eqref{eq:heun1}
has as its solution the Heun polynomial
of order $(n-1)$;
\item the points $\lambda=\lambda_j$ $(j=1,\ldots,n)$ are the branching
points of the hyperelliptic curve $\vG$ \eqref{curve:hyp} for the finite-gap
solution of Heun's equation with characteristics $m_i$ \eqref{eq:mi}.
\end{enumerate}
\end{thm}

\begin{proof}
For
\begin{equation*}
m_0=m_1=m_2=m_3=0
\end{equation*}
the statement of the Theorem is easily checked by the straightforward
substitution of the polynomial of zero order into the Heun equation.

Now, let at least one of the characteristics $m_i$ be different from zero.
The first part of the statement is then a consequence of the general theory
of Heun's polynomials (see, e.g., \cite{HDE, SL}).

The second part will be proved assuming the opposite.
Let for some value $\lambda=\lambda_1$ a solution of Heun's
equation with the characteristics $m_i$ \eqref{eq:mi} be a polynomial
$F(z)$ of order $(n-1)$. Also, assume that the point $\lambda=\lambda_1$
is not a branching point of the hyperelliptic curve $\vG$ \eqref{curve:hyp}
for the finite-gap solution of Heun's equation with these characteristics.
Then one obtains that the 3-rd order equation \eqref{eq:3} with
nonnegative integer characteristics $\wtm_i\geqslant0$
\begin{equation*}
\wtm_i \in\{m_i,-m_i-1\},\qquad i=0,1,2,3,
\end{equation*}
has as its solution two polynomials of zero order:
\begin{align*}
\Pp_1(z)={}&\Psi(\wtl_1,z),\\
\Pp_2(z)={}&z^{M_{1}}(z-1)^{M_{2}}(z-a)^{M_{3}}F^2(z),\\
M_{i}={}&\wtm_{i}-m_i, \qquad i=1,2,3,\\
\wtl_1={}&\lambda_1-2m_1M_{3}-2m_3M_{1}+2M_{1}M_{3}+M_{1}+M_{3}-{}\\
&{}\quad-(2m_1M_{2}+2m_2M_{1}-2M_{1}M_{2}-M_{1}-M_{2})a.
\end{align*}
Moreover, because the point $\lambda=\lambda_1$ of the curve
$\vG$, as was assumed, is not a branching point, the polynomial
$\Pp_1(z)$ does not have multiple zeros and it also does not have
zeros at the singularities of Heun's equation (cf. Corollary
\ref{cor:heun:zero1}). Hence the polynomials $\Pp_1(z)$ and
$\Pp_2(z)$ are linearly independent. It is not difficult to show that
the equation \eqref{eq:3} can not simultaneously have as its solutions
the polynomials $\Pp_1(z)$ and $\Pp_2(z)$. This means that the
point $\lambda=\lambda_1$ is a branching point of the hyperelliptic curve
$\vG$ \eqref{curve:hyp} with the characteristics $m_i$
\eqref{eq:mi}.
\end{proof}

Based on the statements of the Lemma \ref{lem:nk},
the Theorem \ref{thm:polyheun} and Corollary \ref{cor:lj}
one can suggest one more method of finding the branching points
of the curve $\vG$ \eqref{curve:hyp}.

\begin{cor}
In order to find all branching points $\lambda_j$ of the hyperelliptic
curve $\vG$ \eqref{curve:hyp} associated with `finite-gap' solution
of Heun's equation with the characteristics $m_i$ it is necessary
and sufficient to find all Heun's polynomials being the solutions
of Heun equations (from one to four) with the characteristics
$\wtm_i$ \eqref{eq:nk}.
\end{cor}

The algebraic genus $g$ of the curve $\vG$ \eqref{curve:hyp} can be
calculated by finding a sum of all $n_k\in\Bn$ \eqref{eq:nk}, i.e.
by counting all the branching points $\lambda_j$ $(j=1,\ldots,2g+1)$.

\begin{cor} Let $g$ be an algebraic genus of the hyperelliptic curves:
a) $\vG$ \eqref{curve:hyp}, which is associated with the `finite-gap'
solution of Heun's equation, and
b) $\wtG$ \eqref{curve:hyp1}, which is a spectral curve of the
Treibich-Verdier potential. If all characteristics $m_i$ $(i=0,1,2,3)$
are nonnegative integer then:
\begin{enumerate}
\item in the case of even $N$ \eqref{eq:N}
\begin{equation} \label{eq:g1}
g=\max\left\{\max_{0\leqslant i\leqslant3}m_i, \frac N2
-\min_{0\leqslant i\leqslant3}m_i\right\};
\end{equation}
\item in the case of odd $N$
\begin{equation} \label{eq:g2}
g=\max\left\{\max_{0\leqslant i\leqslant3}m_i, \frac
{N+1}2\right\}.
\end{equation}
\end{enumerate}
\end{cor}

\section{Monodromy group of the finite-gap solutions}

Let $a>1$, $\lambda\neq\lambda_j$, $\im\nu(\lambda)=0$,
$m_i\geqslant 0$. In the vicinity of points $z=0$ and $z=1$
let us consider the following `finite-gap' solutions of Heun's
equation:
\begin{subequations}\label{eq:h:def}
\begin{alignat}{1}
\begin{split}
&y_1(z)=\sqrt{\Psi_{g,N}(\lambda,z)}\cos\left(\frac{\nu(\lambda)}2
\int_0^z\frac{t^{m_1}(t-1)^{m_2}(t-a)^{m_3}\,d
t}{\Psi_{g,N}(\lambda,t) \sqrt{t(t-1)(t-a)}}\right),\\
&y_1(z)=y_{10}+O(z),\qquad z\lra 0;
\end{split}\\
\begin{split}
&y_2(z)=\sqrt{\Psi_{g,N}(\lambda,z)}\sin\left(\frac{\nu(\lambda)}2
\int_0^z\frac{t^{m_1}(t-1)^{m_2}(t-a)^{m_3}\,d
t}{\Psi_{g,N}(\lambda,t) \sqrt{t(t-1)(t-a)}}\right),\\
&y_2(z)=z^{m_1+1/2}(y_{20}+O(z)),\qquad z\lra 0;
\end{split}\\
\begin{split}
&y_3(z)=\sqrt{\Psi_{g,N}(\lambda,z)}\cos\left(\frac{\nu(\lambda)}2
\int_1^z\frac{t^{m_1}(t-1)^{m_2}(t-a)^{m_3}\,d
t}{\Psi_{g,N}(\lambda,t) \sqrt{t(t-1)(t-a)}}\right),\\
&y_3(z)=y_{30}+O(1-z),\qquad z\lra 1;
\end{split}\\
\begin{split}
&y_4(z)=\sqrt{\Psi_{g,N}(\lambda,z)}\sin\left(\frac{\nu(\lambda)}2
\int_1^z\frac{t^{m_1}(t-1)^{m_2}(t-a)^{m_3}\,d
t}{\Psi_{g,N}(\lambda,t) \sqrt{t(t-1)(t-a)}}\right),\\
&y_4(z)=(1-z)^{m_2+1/2}(y_{40}+O(1-z)),\qquad z\lra 1;
\end{split}
\end{alignat}
\end{subequations}
where $y_{10}$, $y_{20}$, $y_{30}$ and $y_{40}$ are some functions
of $\lambda$.

From the Corollary \ref{cor:heun:zero1} it follows that the
function $\Psi_{g,N}(\lambda,z)$ does not have zeros in the real
interval $z\in[0,1]$. Therefore, the connection matrix of the
solutions $y_1(z)$, $y_2(z)$ and $y_3(z)$, $y_4(z)$,
\begin{equation*}
\begin{pmatrix} y_1\\ y_2 \end{pmatrix}=T_{01}
\begin{pmatrix} y_3\\ y_4 \end{pmatrix}
\end{equation*}
can be found relatively simple
\begin{equation} \label{eq:T:def}
T_{01}=\begin{pmatrix}
\cos\vphi(\lambda)&
-\sin\vphi(\lambda)\\
\sin\vphi(\lambda)&
\cos\vphi(\lambda)
\end{pmatrix},
\end{equation}
where
\begin{equation} \label{eq:phi:def}
\vphi(\lambda)=\frac{\nu(\lambda)}2
\int_0^1\frac{t^{m_1}(t-1)^{m_2}(t-a)^{m_3}\,d
t}{\Psi_{g,N}(\lambda,t) \sqrt{t(t-1)(t-a)}}.
\end{equation}

From \eqref{eq:h:def}--\eqref{eq:phi:def} it is not difficult to
obtain the expressions for two generators of the monodromy group
for `finite-gap' solutions of Heun's equation \eqref{eq:heun1}
which correspond to simple loops starting at the point $z=0$ and
positively encircling the points $z=0$ and $z=1$, respectively
\begin{align}
&M_0=\begin{pmatrix}1&0\\0&-1\end{pmatrix},\\
&M_1=T_{01}\begin{pmatrix}1&0\\0&-1\end{pmatrix}T_{01}^{-1}
=\begin{pmatrix}
\cos2\vphi(\lambda)&
\sin2\vphi(\lambda)\\
\sin2\vphi(\lambda)&
-\cos2\vphi(\lambda)
\end{pmatrix}.
\end{align}
Analogously, one can find the third generator
\begin{gather}
M_2=T_{02}\begin{pmatrix}1&0\\0&-1\end{pmatrix}T_{02}^{-1}
=\begin{pmatrix} \cos2\psi(\lambda)&
\sin2\psi(\lambda)\\
\sin2\psi(\lambda)& -\cos2\psi(\lambda)
\end{pmatrix},\\
T_{02}=\begin{pmatrix} \cos\psi(\lambda)&
-\sin\psi(\lambda)\\
\sin\psi(\lambda)& \cos\psi(\lambda)
\end{pmatrix},\\
\psi(\lambda)=\frac{\nu(\lambda)}2
\int_{\Cc}\frac{t^{m_1}(t-1)^{m_2}(t-a)^{m_3}\,d
t}{\Psi_{g,N}(\lambda,t) \sqrt{t(t-1)(t-a)}}.
\end{gather}
Here $\Cc$ is a path in the complex plane connecting
the points $z=0$ and $z=a$ and avoiding zeros of the
polynomial $\Psi_{g,N}(\lambda,z)$.

By exactly the same way one can find the generators
of the monodromy group for the `finite-gap' solutions of
Heun's equation with negative values of the
characteristics $m_i$.

\section*{Concluding remarks\protect\footnote{the author learned
about the work \cite{Tak01} from V.B.Kuznetsov}}

Recently, there have appeared papers \cite{Tak00,Tak01} in which,
following \cite{ODE,VTF,WW}, the solutions of Lam\'e
\eqref{eq:lame} and Treibich-Verdier \eqref{eq:shr},
\eqref{pot:tv} equations are studied by making use of periodic
solution of a differential equation for the product of
eigenfunctions of Schr\"odinger's operator with periodic potential
\cite{ODE,VTF,WW}. As seen from these papers, the author is not
familiar with the results of Darboux \cite{Dbk1,Dbk2,D882},
Krichever \cite{Kri80}, Treibich and Verdier
\cite{Ver88,Trb89,TV90,TV90a,TV92}, Belokolos and Enol'skii
\cite{IE86,BBME,BBEIM,BE89b,BE89a,BE94,EK94,EE94,EE95}, Gesztesy
\cite{GW94d,GW94a,GW94b,GW94c,GW94e} who investigated in detail
the finite-gap solutions of these equations as well as their
connection with the Calogero-Moser system with elliptic
interaction \footnote{notice also that the whole volume 36 of {\it
Acta Appl.Math.} for 1994 is dedicated to the memory of Verdier
and contains many works directly or indirectly related to the
Treibich-Verdier equation}.

Despite this, in the part of studying the properties of solutions
of Lam\`e and Treibich-Verdier equations with integer
characteristics $m_i$ the author of those papers not only
re-discovered already known facts, such as the finite number of
eigenvalues for the periodic and anti-periodic problems, but for
two special cases he calculated their number. In this part, the
results of the papers \cite{Tak00,Tak01} coincide with those of
the present paper (cf. Corollary 9), taking into account that the
eigenvalues for periodic and anti-periodic problems are the gap
edges of the spectrum of finite-gap periodic potential (see, for
instance, \cite{Tit2}) or, which is the same, are the branching
points of the hyperelliptic spectral curve $\vG$
\eqref{curve:hyp}. Unfortunately, in \cite{Tak00,Tak01} in the
general case the formula for the number of eigenvalues is an
unproven hypothesis (although the right one as follows from the
results of the present work), and apart from that there are no
methods suggested for calculating those eigenvalues.

Moreover, from our point of view, the analysis of finite-gap
solutions of Lam\`e and Treibich-Verdier equations, in contrast to
the case of Heun's equation, should be carried out not by making
use of the product of eigenfunctions of Schr\"odinger's operator
but by using the method of Krichever-Hermite ansatz
\cite{BBME,BBEIM,BE89b,BE94}. This is because, in the first case a
solution is expressed in terms of elliptic functions in $N$
\eqref{eq:N} unknown parameters solving $N$ transcendental
equations. In the second case, one needs to solve only an
algebraic equation of $N$th order.

Although, perhaps the formulas used in \cite{Tak00,Tak01}
serve best the problem considered there, i.e. to apply the method
of Bethe ansatz to the Treibich-Verdier equation aiming to study
the one-particle Inozemtsev's model.

\appendix
\section{Simplest `finite-gap' solutions}

At the end of this paper we would like to list polynomials
$\Psi(\lambda,z)$ and canonical equations \eqref{curve:hyp}
of the hyperelliptic curves $\vG=\{(\nu,\lambda)\}$ for some
simplest `finite-gap' solutions \eqref{eq:Y} of Heun's
equation with nonnegative characteristics $m_i$ $(i=0,1,2,3)$.
Our examples are indexed by the characteristics
$(m_0,m_1,m_2,m_2)$.

\begin{description}
\item[(0,0,0,0)]
\begin{align*}
&\Psi(\lambda,z)=1,\\ &\nu^2=\lambda.
\end{align*}
\item[(1,0,0,0)]
\begin{align*}
&\Psi(\lambda,z)=\lambda+z-a-1,\\
&\nu^2=(\lambda-1)(\lambda-a)(\lambda-a-1).
\end{align*}
\item[(0,1,0,0)]
\begin{align*}
&\Psi(\lambda,z)=z\lambda+a,\\
&\nu^2=\lambda(\lambda+a)(\lambda+1).
\end{align*}
\item[(1,1,0,0)]
\begin{align*}
&\Psi(\lambda,z)=z\lambda+z^2+a,\\
&\nu^2=(\lambda+a+1)(\lambda^2-4a).
\end{align*}
\item[(0,0,1,1)]
\begin{align*}
&\Psi(\lambda,z)=(z-1)(z-a)\lambda+(a+1)z^2-4az+a(a+1),\\
&\nu^2=(\lambda+a+1)(\lambda^2-4a).
\end{align*}
\item[(1,1,1,0)]
\begin{align*}
&\Psi(\lambda,z)=z(z-1)\lambda^2
+\Big(z^3+3(a-1)z^2-3(a-1)z-a\Big)\lambda
-3a(a-1),\\
&\nu^2=\lambda(\lambda+3a)(\lambda+3(a-1))
(\lambda^2+2(2a-1)\lambda-3).
\end{align*}
\item[(1,1,1,1)]
\begin{align*}
&\Psi(\lambda,z)=z(z-1)(z-a)\lambda+(z^2-a)^2,\\
&\nu^2=\lambda(\lambda+4a)(\lambda+4).
\end{align*}
\item[(2,0,0,0)]
\begin{align*}
&\Psi(\lambda,z)=\lambda^2+(3z-5(a+1))\lambda+9z^2-12(a+1)z
+(a+4)(4a+1),\\
&\nu^2=(\lambda-a-1)(\lambda-a-4)(\lambda-4a-1)
(\lambda^2-4(a+1)\lambda+12a).
\end{align*}
\item[(0,2,0,0)]
\begin{align*}
&\Psi(\lambda,z)=z\lambda^2 +\Big(3z^2-3(a+1)z+a\Big)\lambda
+9z^3-9(a+1)z^2-3a(a+1),\\
&\nu^2=(\lambda+3)(\lambda+3a)(\lambda+3(a+1))
(\lambda^2+4(a+1)\lambda+12a).
\end{align*}
\item[(2,1,0,0)]
\begin{align*}
&\Psi(\lambda,z)=z^2\lambda^2+3\Big((a+1)z+a\Big)z\lambda+9az^2+9a^2,\\
&\nu^2=(\lambda-3(a+1))(\lambda^2-2\lambda-3(4a+1))
(\lambda^2-2a\lambda-3a(a+4)).
\end{align*}
\item[(0,0,2,1)]
\begin{align*}
\begin{split}
{}&\Psi(\lambda,z)=(z-1)^2(z-a)\lambda^2+{}\\
{}&\qquad+\Big(4az^2-(a+3)(3a-1)z+2a(3a-1)\Big)(z-1)\lambda-{}\\
{}&\qquad-12az^3+9(a+1)^2z^2-36az-3a(3a^2-6a-1),
\end{split}\\
{}&\nu^2=(\lambda-3)(\lambda^2+4a\lambda-12a)
(\lambda^2+2(3a-1)\lambda+3(3a^2-6a-1)).
\end{align*}
\item[(2,1,1,0)]
\begin{align*}
\begin{split}
{}&\Psi(\lambda,z)=z(z-1)\lambda^2
+\Big(3z^3+(3a-10)z^2-(3a-8)z-a\Big)\lambda+{}\\
{}&\qquad+9z^4-24z^3-2(3a-8)z^2-a(3a-8),
\end{split}\\
{}&\nu^2=(\lambda+3a-8)(\lambda+3a+1)
(\lambda^3+4(a-2)\lambda^2-16a\lambda+64a).
\end{align*}
\item[(0,1,1,2)]
\begin{align*}
\begin{split}
{}&\Psi(\lambda,z)=z(z-1)(z-a)^2\lambda^2+{}\\
{}&\qquad+\Big(3(2a+3)z^3-(3a^2+18a+8)z^2
+a(3a+10)z+a^2\Big)(z-a)\lambda+{}\\
{}&\qquad+9a(a+3)z^4-24a(3a+1)z^3+18a^2(a+3)z^2-3a^3(a+3),
\end{split}\\
{}&\nu^2=(\lambda+3a)(\lambda+3(a+3))
(\lambda^3+4(a+4)\lambda^2+16(3a+4)\lambda+192a).
\end{align*}
\item[(2,1,1,1)]
\begin{align*}
\begin{split}
{}&\Psi(\lambda,z)=z(z-1)(z-a)\lambda^3+{}\\
{}&\qquad+\Big(3z^4-4(a+1)z^3+2(a+1)^2z^2-2a(a+1)z+a^2\Big)\lambda^2+{}\\
{}&\qquad+\Big(9z^5-15(a+1)z^4-5(a^2-7a+1)z^3
+15(a-1)^2(a+1)z^2-{}\\
{}&\qquad\qquad-15a(a-1)^2z-2a^2(a+1)\Big)\lambda-15a^2(a-1)^2,
\end{split}\\
\begin{split}
{}&\nu^2=\lambda(\lambda^2-2(a+1)\lambda-15(a-1)^2)\times{}\\
{}&\qquad\times(\lambda^2-2(a-2)\lambda-15a^2)
(\lambda^2+2(2a-1)\lambda-15).
\end{split}
\end{align*}
\item[(3,0,0,0)]
\begin{align*}
\begin{split}
{}&\Psi(\lambda,z)=\lambda^3+2\Big(3z-7(a+1)\Big)\lambda^2+{}\\
{}&\qquad+\Big(45z^2-78(a+1)z+49a^2+158a+49\Big)\lambda+{}\\
{}&\qquad+225z^3-405(a+1)z^2+9(3a+8)(8a+3)z-12(a+1)(3a^2+26a+3),
\end{split}\\
\begin{split}
{}&\nu^2=(\lambda-4(a+1))
(\lambda^2-10(a+1)\lambda+3(3a^2+26a+3))\times{}\\
{}&\qquad\times(\lambda^2-2(2a+5)\lambda+3(8a+3))
(\lambda^2-2(5a+2)\lambda+3a(3a+8)).
\end{split}
\end{align*}
\end{description}

\end{document}